\documentclass[10pt,a4paper,final,fleqn,pdftex]{scrartcl} % KOMA-Script article scrartcl
\usepackage{lipsum}
\usepackage{url}
\usepackage[nochapters]{classicthesis} % nochapters

\usepackage[notcite,notref]{showkeys}
\usepackage{amsthm}
\usepackage{amsfonts}
\usepackage{amsmath}
\usepackage{amssymb}
\usepackage{graphicx}
\graphicspath{{./}}
\usepackage{nicefrac}
\usepackage[smaller]{acronym}

\usepackage[english]{babel}
\usepackage[T1]{fontenc}
\usepackage{microtype}
\usepackage[affil-it]{authblk}
\newcounter{thmcount}
\newenvironment{thm}
{% This is the begin code
\refstepcounter{thmcount}{\endgraf\noindent\spacedlowsmallcaps{Theorem}}
\arabic{thmcount} \begin{it}
}
{% This is the end code
\end{it}\smallskip}
\newenvironment{lmm}
{% This is the begin code
\refstepcounter{thmcount}{\endgraf \noindent\spacedlowsmallcaps{Lemma}}
\arabic{thmcount} \begin{it}
}
{% This is the end code
\end{it}\smallskip}
\newenvironment{crl}
{% This is the begin code
\refstepcounter{thmcount}{\endgraf\noindent\spacedlowsmallcaps{Corollary}}
\arabic{thmcount} \begin{it}
}
{% This is the end code
\end{it}\smallskip}

\theoremstyle{remark}
\newtheorem{rmk}{Remark}%[section]
\numberwithin{equation}{section}
\newcommand{\ie}{i.\,e.}

\newcommand{\eg}{e.\,g.}

\newcommand{\prob}[1]{\ensuremath{P\left(#1\right)}}
\newcommand{\trans}[2]{\ensuremath{\mathcal{P}\left(#1, #2\right)}}
\newcommand{\mean}[1]{\ensuremath{\mathbb{E}\left[#1\right]}}
\newcommand{\marginal}{\ensuremath{P_{l}}}
\newcommand{\alphamax}{\ensuremath{\alpha_{\text{\tiny max}}}}
\newcommand{\kmax}{\ensuremath{k_{\text{\tiny max}}}}

\usepackage{mathtools}

\DeclarePairedDelimiter\floor{\lfloor}{\rfloor}

\begin{document}
    % \title{\rmfamily\normalfont\spacedallcaps{the title}}
    % \author{\spacedlowsmallcaps{tyler durden}}
    % \date{} % no date
    \title{\rmfamily\normalfont\spacedallcaps{Asymptotics for the
        Late Arrivals Problem}}

    \author[1]{Carlo Lancia}
    \author[2]{Gianluca Guadagni}
    \author[3]{Sokol Ndreca}
    \author[4]{Benedetto Scoppola}
    \affil[1]{Leiden University}
    \affil[2]{University of Virginia}
    \affil[3]{Universidade Federal de Minas Gerais}
    \affil[4]{Universit\`a di Roma `Tor Vergata'}

    \maketitle

  \begin{abstract}
    \noindent
    We study a discrete time queueing system where deterministic arrivals have
    i.i.d.\ exponential delays $\xi_{i}$.
    The standard deviation $\sigma$ of the delay is finite,
    but much larger than the deterministic unit interarrival time.
    We describe the model as a bivariate Markov chain,
    we prove that it is ergodic and then
    we focus on the unique joint equilibrium distribution.
    We write a
    functional equation for the bivariate generating function,
    finding the solution of such equation on a subset of its
    set of definition. This solution allows us to
    prove that the equilibrium distribution of the Markov chain
    decays super-exponentially
    fast in the quarter plane.
    Finally, exploiting the latter result,
    we discuss the numerical computation of the stationary
    distribution, showing the effectiveness of a simple approximation
    scheme in a wide region of the parameters.
    The model, motivated by air and railway traffic, was
    proposed many decades ago by Kendall~\cite{Kendall1964} with the
    name of ``late arrivals problem'', but no complete solution has been found so far.
  \end{abstract}

    %\tableofcontents

    \section{Introduction}
    \label{sec:intro}
      In this paper we consider a single-server queue with deterministic service
      time, which is assumed of unitary length for the sake of simplicity.
      The $i$th customer arrives to the system at time
      \begin{equation}
      \label{eq:eda1}
      t_i=i+\xi_i\,, \qquad i\in\mathbb{N}\,,
      \end{equation}
      where  $\{\xi_i\}$ are i.i.d. exponential random variables with parameter $\beta$.

      In the limit  $\beta\to0$ the point process~(\ref{eq:eda1})
      weakly converges to a
      Poisson process of parameter $1$,
      whereas for fixed $\beta$ the arrivals are
      negatively autocorrelated, see~\cite{collings1976limit, gns} and
      references therein.
      \begin{rmk}
      Although the results in~\cite{gns} are stated under the hypothesis that
      the probability density function of the delays $\{\xi_i\}$ has
      compact support, this assumption does not play
      any role in establishing the convergence to a Poisson process.
      As a consequence, the very same result applies here too.
      \end{rmk}

      We study the system described above for fixed $\beta$ and
      we assume that arriving customers might balk with independent probability
      $1-\rho$.
      In other words, each customer can be
      deleted independently with probability $1-\rho$ before
      joining the queue.
      Besides being a mathematical expedient that ensures
      the existence of a stationary state\footnote{See
      Lemma~\ref{rmk:uniqueness} below.},
      the balking is mainly a way to model empty intervals
      in a constant stream of customers.
      Again, the point process~(\ref{eq:eda1}) with balking
      weakly converges to a Poisson process, but with parameter
      $\rho$~\cite{gns}.
      In what follows, we refer to the balked version of~(\ref{eq:eda1}) as Exponentially Delayed Arrivals (EDA).
      % After Kendall (see below for his exact words) we name
      % Exponentially Delayed Arrivals (EDA) the thinned version of the
      % arrival process~(\ref{eq:eda1}).

      Service can be delivered by the unique server only at discrete times.
      The length of the queue at time $t$ is $n_t$; it represents the number of
      customers waiting to be served, including  the customer that will be
      served precisely at time $t$, if any.
      Due to the balking, it is immediate to see that the traffic
      intensity of the system
      is given by $\rho$; see~\cite{gns} for details.
      Using Kendall's notation we hereafter refer to the queue model
      described so far as $EDA/D/1$.

      The $EDA/D/1$ model is motivated by the description
      of public and private transportation systems,
      including buses, trains, aircraft~\cite{ball2001analysis,gns,
      gwiggner2010analysis,iovanellaimpact}
      and vessels~\cite{govier1963stock,jagerman2003vessel},
      appointment scheduling in outpatient
      services~\cite{bailey1952study,cayirli2003outpatient,mercer1960,
      mercer1973queues}
      crane handling in dock operations~\cite{daganzo1990productivity,
      edmond1975operation},
      and in general any system where scheduled arrivals are
      intrinsically subject to random variations.
      Preliminary results show that the model described above fits very well
      with actual data of inbound air traffic over a large hub,
      see~\cite{cills,ll2017}.

      The appearance of the stochastic point process~\eqref{eq:eda1}
      can be traced back to Winsten's seminal paper~\cite{win}.
      Winsten named such a queueing model
      \textit{late process} and obtained results for
      the special case $\xi_i\in [0, 2]$ and service
      time exponentially distributed.
      At the end of~\cite{win} there is a discussion on Winsten's
      results by Lindley, Wishart, Foster, and Tak\'acs, where they state
      that Winsten's paper can be considered as the first
      treatment of a queueing model with correlated
      arrivals~\cite[pages 22-28]{win}.

      The same problem was investigated also by Kendall.
      In~\cite[page~11]{Kendall1964}
      he remarked the great importance of
      systems with arrivals like~\eqref{eq:eda1}:
      ``\textit{[...] perhaps too much attention has been paid
      to rather uninteresting
      variations on the fundamental Poisson stream.
      As soon as one considers
      variations dictated by the exigencies of the real world, rather than
      by the pursuit of mathematical elegance,
      severe difficulties are encountered;
      this is particularly well illustrated by the
      notoriously difficult problem
      of \emph{late arrivals}.}''
      Kendall also provided the following elegant interpretation: if the
      random variables $\xi_i$ are non-negative then the
      process defined by~\eqref{eq:eda1}
      is the output of the stationary $D/G/\infty$ queueing
      system~\cite{Kendall1964}.
      In particular, if the random variables $\xi_i$
      are exponentially distributed then $EDA/D/1$ can be viewed as a
      2-stage tandem queueing network
      \[D/M/\infty \rightarrow \cdot/D/1 \,.\]
      However, this is not the approach followed in this work.

      Some years later, under the hypothesis that
      $\xi_i>0$, Nelsen and Williams
      exactly characterised in~\cite{nw} the distribution of the
      inter-arrival time intervals and their correlations.
      They also gave an explicit expression of these
      quantities in the particular case of
      $\xi_i$'s exponentially distributed.

      After the '70s only approximations of the arrival
      process~\eqref{eq:eda1}
      \cite{bloomfield1972low,sabria1989approximate}
      or numerical studies of its output
      \cite{almaz2012simulation,ball2001analysis,nikoleris2012queueing}
      seem to have appeared in the literature.
      In particular, in~\cite{gns} the authors presented a self-contained
      study of an arrival process like~\eqref{eq:eda1},
      assuming for $\xi_i$
      a compact-support distribution. They also proposed an
      approximation scheme that keeps the correlation of the arrivals and
      is able to compute in a quite accurate way the
      quantitative features of the queue.
      To the best of our knowledge,
      a queueing system with arrivals described by~\eqref{eq:eda1}
      still remains an open problem and the
      best results obtained so far are due to
      Winstein in 1959~\cite{win}.

      \medskip

      $EDA/D/1$ is an example of a queueing system with
      correlated arrivals,
      a subject broadly studied in past years.
      There are many ways to impose a correlation to the arrival process.
      For instance, the parameters of the process may depend on their past
      realisation, as in~\cite{drezner3queue},
      or on some on/off sources, as in~\cite{wittevrongel1999discrete}.
      Another relevant example of a queue model with correlated arrivals
      is the so-called Markov Modulated Queueing System.
      In Markov Modulated Queueing Systems the parameters are driven by
      an independent external Markovian process,
      see~\cite{adan2003single,asmussen2000multi,combe1998bmap,
      lucantoni1991new,pacheco2009markov}
      and references therein.
      Our model shares with Markov Modulated Queueing Systems the property
      that one can define an external and independent
      Markovian process that drives the arrival rates. However,
      we see in Section~\ref{sec:bivgf} that
      the output of this external drive
      has also a deterministic effect on the queue length.
      More precisely, $EDA$ can be interpreted as an independent drawing from an external pool of customers \emph{late} at time $t$,
      see~\eqref{eq:edal_t} below.
      Due to the memoryless property of the exponential delays,
      each customer late at time $t$ will be still in the pool
      at time $t+1$ independently with probability $q\equiv e^{-\beta}$.
      This leads to binomial transitions in the number of late customers.

      In Section~\ref{sec:bivgf} we show that $EDA/D/1$ can be
      described as a bivariate Markov chain
      representing the queue length and the number of late
      customers. We prove that such a bivariate chain is ergodic and
      write the balance equations of the stationary distribution,
      finding a functional equation for the bivariate
      generating function.

      There exists an extensive literature about two-dimensional
      Markov models. Many methods for attacking the problem are available
      under two assumptions, namely,
      spatial homogeneity and finiteness of at least one marginal chain,
      see~\cite{bini2005numerical,
      gail2000use,grassmann2002real,
      latouche1999introduction,mitrani1995spectral,
      neuts1989structured,neuts1995matrix}.
      Unfortunately, the Markov chain defined in
      Section~\ref{sec:bivgf} does not satisfy any of the mentioned
      requirements.

      When both components of the Markov chain are infinite but
      space homogeneity is still ensured,
      the problem is typically attacked by reduction to a
      Riemann-Hilbert Boundary Value Problem.
      These may be solved, for example, by
      the uniformisation technique \cite{kingman1961two},
      conformal mappings~\cite{cohen1988boundary,
      cohen1983boundary,fayolle1979two},
      the compensation meth\-od~\cite{adan1993compensation},
      or the Power Series Approximation \cite{blanc1987numcoupled,
      blanc1987numerical,koole1997use,hooghiemstra1988power}.
      % Usually, Power Series Approximation is used to obtain the generating function in
      % terms of a power series in the load $\rho$ although
      % different parameters may be used, see~\cite{walraevens2010power}.
      % A power series expansion in a parameter different
      % from $\rho$ is also the strategy we adopt in what follows.
      The aforesaid binomial transitions are responsible for the lack of spatial homogeneity and are often encountered in Mathematical Biology~\cite{artalejo2007evaluating,
      brockwell1982birth,economou2004compound}.

      To the best of our knowledge, the functional equation~\eqref{eq:eda4} below
      has never previously appeared in the literature.
      Yet it is possible to mark some analogies with the functional equation
      in~\cite{economou2009q,economou2010single,kapodistria2011m},
      the most important being that in both equations the right hand side
      exhibits the generating function computed in a convex
      combination in the parameter $q=e^{-\beta}$.
      Other examples of chains with binomial transitions
      may be found
      in~\cite{adan2009synchronized,altman2006analysis,
      neuts1994interesting,yechiali2007queues}.

      In Section~\ref{sec:marg-distr-late} we study the marginal
      distribution of the number of late customers and we obtain its exact
      analytical expression, which reveals the rich
      combinatorial structure of the problem.
      This intermediate result allows us to
      show that the stationary distribution of the $EDA/D/1$
      queue has a super-exponential decay.
      Finally, in Section~\ref{sec:numerics} we show that such a
      super-exponential behaviour enables a simple, yet very effective,
      numerical approximation scheme of the system balance equations.
      For a wide range of the system parameters, including typical values
      for real traffic applications of the model, we give a very good
      a priori estimate of the total-variation distance between the true
      and the approximate solution.

    \section{Stationary distribution:
    generating function and balance
    equations}
    \label{sec:bivgf}

      Let us consider the process $n_t$, which describes the length of
      the queue at time $t$. This process is governed by
      the stochastic recursion
      \begin{equation}
      \label{eq:eda2a}
      n_{t+1}= n_t+m_{(t,t+1]} -(1-\delta_{n_t,0}) \,,
      \end{equation}
      where $m_{(t,t+1]}$ is the number of arrivals in the interval
      $(t,t+1]$ according to the arrival process~\eqref{eq:eda1}, and
      $\delta_{i,j}$ is the usual Kronecker's delta.
      The term $1-\delta_{n_t,0}$ represents the action of the server in
      decreasing the queue length by one if at time $t$ the queue is non-empty.
      Since the service time is deterministic, we focus on the so-called \emph{embedded process} by observing the system at $t\in\mathbb{N}$, \ie{} at departure instants.

      The quantity $m_t \equiv m_{(t,t+1]}$ depends in general on the
      whole previous history of the system.
      Indeed, if for some large value of $T$,
      $m_{s}=0$ for any $s\in \{t-T, t-T+1,...,t-1\}$
      then $m_{t}$ is large with great probability.
      Conversely, if in the recent
      past the values of $m_{s}$ have been large then $m_{t}$
      is expected to be small.
      This suggests that the arrival process
      is negatively autocorrelated, as proven in~\cite{gns}.
      Hence, the recursion~\eqref{eq:eda2a}
      does not depend only on the present value of $n_t$, and
      the memory of the process is infinite
      since $T$ can be arbitrarily large.

      Let us now denote by $l_t$ the number
      of customers that are late at time $t$, that is
      to say,
      \begin{equation}
      \label{eq:edal_t}
      l_t \equiv \big\vert\{ 0 \leq i \leq t \text{ such that }
      \xi_i > t-i \}\big\vert \,.
      \end{equation}
      Let us next define $p \equiv \int_0^1f_\xi(t)dt=\int_0^1\beta e^{-\beta t}dt=1-e^{-\beta}$ and $q \equiv e^{-\beta}$.
      Given the value of $l_t$, the random variable
      $m_t$ is binomial with parameters $l_t$ and $p$.
      According to the memoryless property of the exponential delays $\{\xi_i\}$,
      the number of unit time intervals in which a customer is late is a geometric random variable with parameter $q$.
      In other words, a customer will be late for $k$ consecutive time intervals with probability $q^kp$.
      Hence, the embedded process $l_t$, $t \in \mathbb{N}$, is a discrete-time Markov chain.
      If the customer scheduled in the interval $(t,t+1]$ has balked then
      \begin{equation}
      \label{eq:edabinom1}
      \prob{m_t=j \,|\, l_t = l}= \binom{l}{j} p^j q^{l-j} \equiv b_{j,l} \,,
      \end{equation}
      otherwise
      \begin{equation*}
      %\label{eq:edabinom2}
      \prob{m_t=j  \,|\, l_t = l}= \binom{l+1}{j} p^j q^{l+1-j} = b_{j,l+1}\,.
      \end{equation*}
      All in all,
      \begin{equation}
      \prob{m_t=j \,|\, l_t = l}
      = b_{j,l}\,(1-\rho) + b_{j,l+1}\, \rho \,. \label{eq:edabinom22}
      \end{equation}

\begin{figure}[tbp]
      \centering
      \includegraphics[width=.85\textwidth]{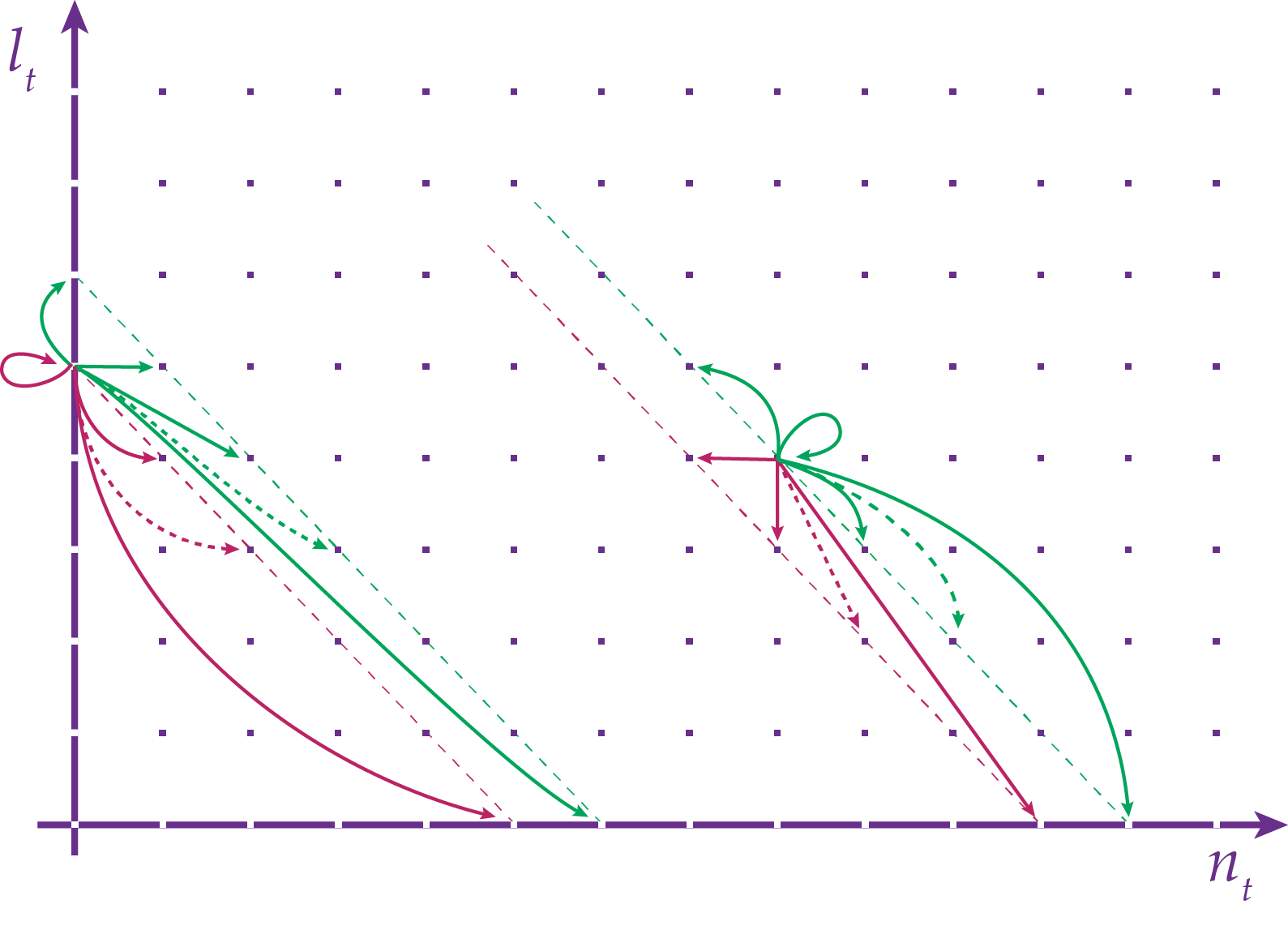}
      \caption[Transitions of $EDA/D/1$ in the quarter plane.]{Transitions of the $EDA/D/1$
      queueing system in the quarter plane.  They happen
      along the lines of Cartesian equation $x+y = n_t+l_t$ and $x+y =
      n_t+l_t - 1$ if $n_t\ne0$, and along the lines of Cartesian equation $x+y = n_t+l_t$ and $x+y =
      n_t+l_t + 1$ if $n_t=0$. Green transitions happen with probability $\rho$
      (no balking), red transitions happen with probability $1-\rho$
      (balking).}
      \label{fig:quarterplane}
      \end{figure}

      We describe the $EDA/D/1$ queue by the embedded process $(n_t, l_t)$, $t\in\mathbb{N}$, a discrete-time Markov chain with the following transition probabilities:
      \begin{align}
      &\text{For } n>0, \label{eq:eda120}\\
      &\quad \trans{(n,l)}{(n+a-1,l-a+1)}= \rho\,b_{a,l+1}\,, &\;\; &0\leq a\leq l\!+\!1\,, \nonumber\\
      &\quad \trans{(n,l)}{(n+a-1,l-a)} = (1-\rho)\,b_{a,l}\,, &\;\; & 0\leq a\leq l\,; \nonumber\\
      &\text{For } n=0, \label{eq:eda121}\\
      &\quad \trans{(0,l)}{(a,l-a+1)}= \rho\,b_{a,l+1}\,, &\;\; &0\leq a\leq l\!+\!1\,, \nonumber\\
      &\quad \trans{(0,l)}{(a,l-a)} = (1-\rho)\,b_{a,l}\,, &\;\; & 0\leq
      a\leq l\,. \nonumber
      \end{align}
      Figure~\ref{fig:quarterplane} displays the transitions of the embedded chain $(n_t, l_t)$ in the quarter plane.

      We now address the existence of the stationary state of the embedded chains $l_t$ and $(n_t, l_t)$.
      \begin{lmm}
      The chain $l_t$ is ergodic if and only if $q<1$.
      \end{lmm}
      \begin{proof}
      If $q<1$, the chain $l_t$ is clearly irreducible. In order to
      prove its positive recurrence,
      we use Foster's criterion \cite[Cor.\ 8.7]{robert2003stochastic}
      setting $f(l_t) = l_t + 1$ as the Lyapunov function.
      Thus, we need to show that there exist suitable positive constants
      $K, \gamma$ such that
      \begin{enumerate}
      \item $\mean{f(l_1) - f(l_0) \,|\, l_0=L} \leq -\gamma$ for $f(L) > K$;
      \item $\mean{f(l_1) \,|\, l_0=L} < \infty$ for $f(L) \leq K$;
      \item the set $\{l \geq 0 \,:\, f(l) \leq K\}$ is finite.
      \end{enumerate}
      First, by~\eqref{eq:edabinom22},
      \begin{align}
      \mean{f(l_1) - f(l_0) \,|\, l_0=L} &= Lq(1-\rho) + (L+1)q\rho - L\,,\nonumber\\
      &= \rho q - L(1-q)\,. \label{eq:foster2}
      \end{align}
      Therefore, point 1 is satisfied, for instance, by the choice $\gamma =1$ and $K=1+\frac{1+\rho q}{1-q}$.
      Next, Figure~\ref{fig:quarterplane} shows that at each iteration
      $l_t$ increases at most by one unit and point 2 is satisfied:
      \begin{equation}
      \label{eq:foster1}
      \mean{f(l_1) \,|\, l_0=L} \leq L+2 \,.
      \end{equation}
      By definition of $f(l_t)$ point 3 is also fulfilled.
      On the other hand, for $q=1$ the chain $l_t$ is not ergodic because it is no longer
      irreducible\footnote{For $q=1$, the chain satisfies in fact $l_{t+1} \geq l_t$.}.
      \end{proof}

      \begin{lmm}
      \label{rmk:uniqueness}
      The bivariate chain $(n_t, l_t)$ is ergodic if and only if $q<1$ and $\rho<1$.
      \end{lmm}
      \begin{proof}
      If $q<1$ and $\rho<1$, the bivariate chain $(n_t, l_t)$ is irreducible, see Figure~\ref{fig:quarterplane}.
      Let us consider the process $\alpha_t = n_t + l_t$, which
      represent the diagonal in the quarter
      plane where the point $(n_t, l_t)$ lies on, see Figure~\ref{fig:quarterplane}.
      The bivariate chain has the property that $|\alpha_{t+1}-\alpha_t|\leq1$.
      Equations~\eqref{eq:eda120}-\eqref{eq:eda121} yield
      \begin{align}
      &P(\alpha_{t+1} = \alpha_t + 1 \,|\, n_t\ne0) = 0\,, \label{eq:16}\\
      &P(\alpha_{t+1} = \alpha_t - 1 \,|\, n_t\ne0) = (1-\rho)\,, \label{eq:17}\\
      &P(\alpha_{t+1} = \alpha_t + 1 \,|\, n_t=0) = \rho\,, \label{eq:18}\\
      &P(\alpha_{t+1} = \alpha_t - 1 \,|\, n_t=0) = 0\,.\label{eq:19}
      \end{align}

      In order to prove the positive recurrence of $(n_t, l_t)$,
      we use again Foster's criterion
      setting $f(n_t, l_t) = M\alpha_t + l_t + 1$ with $M=\nicefrac{2}{(1-\rho)}$.
      From~\eqref{eq:foster2} and~\eqref{eq:16}-\eqref{eq:19},
      \begin{align*}
      &\mean{f(n_1,l_1) - f(n_0,l_0) \,|\, n_0=N, l_0=L} \\%\nonumber\\
      &\qquad= M\mean{\alpha_1-\alpha_0 \,|\, n_0=N, l_0=L} +
      \mean{l_1 - l_0 \,|\, l_0=L}\,, \\%\nonumber\\
      &\qquad\leq M\rho\delta_{N,0} - M(1-\rho)(1-\delta_{N,0}) + \rho
      q - L(1-q) \,. %\label{eq:foster4}
      \end{align*}
      Then, using a little algebra, it can be shown that the first point of Foster's criterion is satisfied
      by $\gamma = 1$ and $K > 1 + \frac{(M+1)(1+\rho q + M\rho)}{1-\rho}$ (for example, $K=\frac{7}{(1-\rho)^2(1-q)}$).
      Point 2 of Foster's criterion holds because
      \begin{equation*}
      %\label{eq:foster3}
      \mean{f(n_1,l_1) \,|\, n_0=N, l_0=L} \leq M(N+L+1) + L + 2,
      \end{equation*}
      where we have used the property that $\alpha_t$ has only nearest-neighbour
      transitions and equation~\eqref{eq:foster1}.
      Point 3 is fulfilled by simply considering the definition of $f(n_t, l_t)$.
      For $\rho=1$, the bivariate chain $(n_t, l_t)$
      is not ergodic because it is no longer
      irreducible\footnote{For $\rho=1$, the process $\alpha_t$
      satisfies in fact $\alpha_{t+1} \geq \alpha_t$.}.
      \end{proof}

      % \begin{rmk}
      %   \label{rmk:uniqueness}
      %   Provided $\rho < 1$, the bivariate Markov chain $(n_t, l_t)$
      %   is irreducible and positive recurrent.
      %   The latter property can be easily checked through the so-called
      %   Foster-Lyapunov criteria.
      %   Let us consider the Lyapunov function $V((n_t, l_t)) = l_t +1$.
      %   It is immediate to prove that
      %   \begin{align*}
      %     \mean{V((n_1, l_1))-V((n_0, l_0)) \;|\; l_0=l} = \rho\,q - l(1-q)\,.
      %   \end{align*}
      %   Then, $\mean{V((n_1, l_1))-V((n_0, l_0)) \;|\; l_0=l}$ is easily bounded from above
      %   \begin{itemize}
      %   \item by a positive constant uniformly in $l$;
      %   \item by a negative constant asymptotically for $l\to\infty$.
      %   \end{itemize}
      %   This is sufficient to infer positive recurrence in the whole
      %   quarter plane, see~\cite{tweedie1976criteria}.
      %   The aperiodicity is trivial.
      %   Thus, $(n_t, l_t)$ is an ergodic Markov chain, and there
      %   exists a unique stationary measure $P_{n,l}$.
      %   We note that the irreducibility no longer holds for
      %   $\rho=1$ because in this case the system is unstable and will eventually
      %   leave any finite region of the quarter plane.
      % \end{rmk}

      For $\rho, q < 1$, Lemma~\ref{rmk:uniqueness} guarantees both
      the existence and the uniqueness of
      the stationary distribution
      \begin{equation}
      \label{eq:pnldef}
      P_{n,l} \equiv \lim_{t\to\infty} P(n_t=n, l_t=l) \,.
      \end{equation}
      Let us consider the following bivariate generating function:
      \begin{equation}
      \label{eq:edagf}
      P(z,y) = \sum_{n,l \geq 0} \, z^n \, y^l \, P_{n,l}\,, \quad |z|,|y|
      \leq 1 \,.
      \end{equation}
      We are now ready to prove the main result of this section.

      \begin{thm}
      \label{thm:eda-functeq}
      The bivariate generating function~\eqref{eq:edagf} satisfies
      \begin{align}
      P(z,y) &= \frac{1+\rho\,(\upsilon-1)}{z}\left[(z-1)
      \,P(0,\upsilon)+P(z,\upsilon)\right] \,,\label{eq:eda4}
      \end{align}
      where $\upsilon = \upsilon(z,y) = z + q\,(y-z)$.
      \end{thm}

      \begin{rmk}
      The functional equation~\eqref{eq:eda4} does not admit
      simple or immediate solutions. It is radically different from
      the functional equations typically studied in the
      literature~\cite{cohen1983boundary,fayolle1999random} and it is
      rather special in this respect. A simple solution can be found only in the
      particular case $z=1$, see Section~\ref{sec:marg-distr-late} below.
      \end{rmk}

      \begin{proof}[Proof of Theorem~\ref{thm:eda-functeq}]
      For each $n,l \geq 0$, the balance equations of $EDA/D/1$ are the following:
      \begin{multline}
      P_{n,l} = (1-\rho) \left( \sum_{j=0}^n P_{j+1, l+n-j} \,
      b_{n-j, l+n-j}
      + P_{0, l+n} \, b_{n, l+n}\right) \\
      + \rho\left( \sum_{j=0}^n\! P_{j+1, l+n-j-1} \, b_{n-j,
      l+n-j} + P_{0, l+n-1} \, b_{n, l+n}\right),\label{eq:eda95}
      \end{multline}
      where $b_{j,l}$ are given by~\eqref{eq:edabinom1}
      and we agree that $P_{n,l}=0$ whenever $n,l<0$.
      The special cases $n=0$ and $n=l=0$ respectively lead to
      \begin{align}
      &P_{0,l} = (1-\rho)\left(P_{1,l}+P_{0,l}\right)\,b_{0,l} +
      \rho\left(P_{1,l-1}+P_{0,l-1}\right)\,b_{0,l} \,,\label{eq:eda96}\\
      &P_{0,0} = (1-\rho)(P_{1,0}+P_{0,0}) \,.\label{eq:eda97}
      \end{align}
      To show that~\eqref{eq:eda95}--\eqref{eq:eda97} hold,
      it suffices to write $P_{n_{t+1},l_{t+1}}$ in terms of
      $P_{n_{t},l_{t}}$ and then neglect the time dependency.
      Take for example~\eqref{eq:eda97}:
      the system is found at time $t+1$ in state $(0,0)$, \ie{}
      with empty queue and no late customers, only
      if at time $t$ it was either in state $(0,0)$ or in state $(1,0)$,
      and the $(t+1)$th scheduled customer\footnote{Cf.\ formula~\eqref{eq:eda1}.} is deleted by thinning.
      Indeed, if at time $t$ the system was in state $(0,0)$ then nothing happens
      and the state remains unchanged, whereas if it was
      in state $(1,0)$ then the
      customer in queue is served and at time $t+1$ the
      system is in state $(0,0)$.
      Similarly, there are four cases such that the system is found
      at time $t+1$ in state $(0,l)$,
      \ie{} with an empty queue and $l$ customers late.
      In the first two cases the system is in state $(1,l)$ or in state $(0,l)$ at time $t$,
      the $(t+1)$th customer is deleted, and no one of
      the $l$ late customers arrives in the interval $[t,t+1)$
      (this event has in fact probability $q^l = b_{0,l}$).
      In the remaining cases
      the system is in state $(1,l-1)$ or in state $(0,l-1)$ at time $t$,
      the $(t+1)$th customer is not deleted\footnote{Therefore, the
      $(t+1)$th customer is added
      to the set of the $l-1$ customers that are already late.}, and no one of
      the $(l-1)+1$ late customers arrive in the interval $(t,t+1]$.
      The latter argument gives~\eqref{eq:eda96} while an easy generalisation
      to the case $n\geq1$ leads to~\eqref{eq:eda95}.

      \noindent
      Let us take~\eqref{eq:eda95}, multiply both sides by $z^n\,y^l$, and then
      sum over $n$ and $l$. The summation of all terms multiplied by $(1-\rho)$ yields
      \begin{multline*}
      (1-\rho)\Bigg\{ \sum_{n,l\geq0} \Bigg[ \sum_{j=0}^n
      P_{j+1,l+n-j} {n+l-j \choose n-j} z^j
      (zp)^{n-j}(yq)^{l} \\
      +P_{0,l+n} {l+n \choose n} (zp)^n (yq)^l \Bigg]
      \Bigg\}\,,%\label{eq:95}
      \end{multline*}
      or equivalently,
      \begin{multline*}
      (1-\rho)\Bigg\{ \sum_{j\geq0}\sum_{n\geq j}
      \Bigg[\sum_{l\geq0} P_{j+1,l+n-j} {n+l-j \choose n-j} z^j
      (zp)^{n-j}(yq)^{l} \\
      +P_{0,l+n} {l+n \choose n} (zp)^n (yq)^l \Bigg]
      \Bigg\}\,.%\label{eq:96}
      \end{multline*}
      The change of indices $k=n-j$ and $m=l+n-j=l+k$ yields
      \begin{align*}
      &(1-\rho)\Bigg\{ \sum_{j\geq0}\sum_{n\geq j}
      \Bigg[\sum_{l\geq0} P_{j+1,l+n-j} {n+l-j \choose n-j} z^j
      (zp)^{n-j}(yq)^{l} \Bigg\}\\
      &=(1-\rho)\Bigg\{ \sum_{j\geq0} z^j \sum_{m\geq0} \Bigg[
      \sum_{k=0}^n P_{j+1,m} {m \choose k}
      (zp)^{k}(yq)^{m} \Bigg]\Bigg\}\,,\\%\label{eq:97}\\
      &=\frac{1-\rho}{z} \sum_{j\geq1} \sum_{m\geq0} P_{j,m} z^j
      (zp+yq)^m , \\%\label{eq:98}\\
      &= \frac{1-\rho}{z} \left[ P(z,zp+yq) - P(0,zp+yq)
      \right] \,, %\label{eq:99}
      \end{align*}
      to which we still have to sum the contribution
      \begin{equation*}
      %\label{eq:102}
      (1-\rho)\sum_{n,l\geq0} P_{0,l+n} {l+n \choose n} (zp)^n
      (yq)^l = (1-\rho)P(0,zp+yq) \,.
      \end{equation*}
      All in all, the sum of all terms multiplied by $(1-\rho)$ is
      \begin{align*}
      (1-\rho) \left[ P(0,zp+yq) + \frac{1}{z}\left( P(z,zp+yq) -
      P(0,zp+yq) \right)\right].
      \end{align*}

      \noindent
      In a completely analogous way we can compute the sum of the terms multiplied by $\rho$,
      which turns out to be
      \begin{align*}
      \rho(zp+yq)\Big[ P(0,zp+yq) + \frac{1}{z}\left( P(z,zp+yq) -
      P(0,zp+yq) \right)\Big].
      \end{align*}
      Summing up the two contributions, we get~\eqref{eq:eda4}.
      % \qed
      \end{proof}

      \begin{rmk}
      We end this section with a discussion of the special case $q=0$.
      In this regime, the right-hand side of equation~\eqref{eq:eda4} does not
      depend on $y$ anymore, and $P(z,y) \equiv Q(z)$.  The number of
      late customers is in fact $l_t = 0$ as the $i$th customer can
      not have a delay $\xi_i \geq 1$.
      % The system reduces to a simple $D/D/1$ queue
      % with balking, and the stationary probability to
      % have a void queue is $1-\rho$.
      Then, equation~\eqref{eq:eda4} yields directly
      \begin{equation}
      \label{eq:125}
      Q(z) = \frac{1+\rho(z-1)}{z} [(z-1)Q(0) + Q(z)] \,,
      \end{equation}
      where $Q(0)=1-\rho$ is the stationary probability of a void
      queue. Therefore, equation~\eqref{eq:125} is equivalent to
      \begin{equation*}
      % \label{eq:124}
      Q(z) = 1+\rho(z-1)\,,
      \end{equation*}
      which is the classical result of a $D/D/1$ queue with
      balking.
      \end{rmk}

    \section{The marginal distribution of late customers}
    \label{sec:marg-distr-late}

      In this section we focus on \marginal{}, the marginal distribution of late
      customers. First, we iterate the functional equation~\eqref{eq:eda4}
      to obtain the generating function of \marginal{} in the form of an
      infinite product. Then, we invert the generating function
      and find the exact analytical expression of \marginal{}. Finally, we
      derive the asymptotic behaviour of \marginal{} and use it to infer
      asymptotics for $P_{n,l}$.

      \smallskip

      \noindent
      The marginal distribution of late customers is
      \[\marginal = \sum_{n\geq0} P_{n,l}\]
      and its generating function is
      \[\sum_{l\geq0} \marginal y^l = \sum_{n,l\geq0}  P_{n,l} \, y^l = P(1,y) \,.\]
      Setting $z=1$ into equation~\eqref{eq:eda4} yields
      \begin{equation}
      \label{eq:marginal1}
      P(1,y) = [1+\rho q(y-1)] P(1,1+q(y-1)) \,.
      \end{equation}
      Evaluating the last equation in $(1, 1+q(y-1))$ yields
      \[P(1, 1+q(y-1)) = [1+ \rho q^2(y-1)] P(1,1+q^2(y-1))\,.\]
      Iterating~\eqref{eq:marginal1} $N$ times,
      \begin{equation*}
      %\label{eq:marginal2}
      P(1,y) = \left[ \prod_{k=0}^{N-1} 1+\rho q^{k+1}(y-1) \right] P(1, 1+q^N(y-1))\,.
      \end{equation*}
      The limit of $\prod_{k = 0}^{N-1} [1+\rho q^{k+1}(y-1)]$ for $N\to\infty$
      exists for each $q<1$ and $y \in \mathbb{C}$, see~\cite{ahlfors1953complex}.
      Moreover, $\lim_{N\to\infty}P(1, 1+q^N(y-1))=1$
      Therefore, we have proven the following
      \begin{crl}
      For $q<1$ and $|y|\leq 1$,
      \begin{equation}
      \label{eq:marginal3}
      P(1,y) = \prod_{k \geq 0} \left[ 1+\rho q^{k+1}(y-1) \right] \,.
      \end{equation}
      \end{crl}

      \begin{rmk}
      The infinite product~\eqref{eq:marginal3} has an interesting
      combinatorial interpretation:
      \begin{equation}
      \label{eq:23}
      P(1,y) = \prod_{k \geq 0} \left[ 1+\rho q^{k+1}(y-1) \right] =
      \frac{(\rho (1-y); q)_\infty}{1+\rho(y-1)}\,,
      \end{equation}
      where $(a; q)_\infty = \prod_{k \geq 0} [1 - aq^k]$
      is the \textit{infinite $q$-Pochhammer
      symbol}, also known as \emph{infinite $q$-ascending factorial in $a$}.
      % of the pair $(\rho(1-y),\,q)$.  It is the $q$-analog of the
      % descending factorial, also known as Pochhammer symbol.
      For $y= 1-\nicefrac{q}{\rho}$,
      $P(1,y) = \phi(q) \,(1-q)^{-1}$, where $\phi(q)$ is the well-known {\it Euler function}.
      \end{rmk}
      \begin{rmk}
      \label{rmk:insight}
      For $q<1$, the right-hand side of~\eqref{eq:23} is analytic for each
      $y\in\mathbb{C}$. Therefore, the power series $P(1,y) =
      \sum_{l\geq0} P_l \, y^l$, convergent for each $|y|\leq1$,
      can be analytically continued in the whole complex plane.
      As such, we expect the marginal distribution $P_l$
      to decrease super-exponentially fast in $l$.
      \end{rmk}

      Following the insight given by Remark~\ref{rmk:insight}, we shift now
      the focus to the asymptotic behaviour of $\marginal$ and $P_{n,l}$.
      Expanding the product and rearranging it in powers of $\rho(y-1)$ yields
      \begin{align}
      P(1,y) &= \prod_{k\geq0}\left(1 + \rho \, q^{k+1}(y-1)\right) \nonumber\\
      &= 1 + \sum_{k\geq1} \rho^k (y-1)^k \left[\sum_{m\geq{k+1 \choose 2}}
      d(m;k) q^m \right],
      \label{eq:129}
      \end{align}
      where $d(m;k)$ is the number of partitions of $m$ in $k$ distinct parts.

      The following theorem holds:
      \begin{thm}
      \label{lmm:yl-P1y}
      Let \marginal{} be the equilibrium marginal distribution
      of the number of late customers and $P(1,y)$ its generating function.
      Then,
      \begin{align}
      P(1,y) &= \sum_{k\geq0} \frac{\rho^k q^{{k+1 \choose 2}} (y-1)^k}{\prod_{i=1}^k [1-q^i] }\,, \label{eq:20} \\
      \marginal %&= [y^l]P(1,y) =
      % \frac{1}{l!}\frac{d^l}{dy^l}P(1,y)\Big\vert_{y=0} \nonumber\\
      &= \sum_{k\geq l} \frac{(-1)^{k-l}\rho^k q^{{k+1 \choose 2}}
      {k\choose l}}{\prod_{i=1}^k [1-q^i] } \,. \label{eq:138}
      \end{align}
      \end{thm}

      Theorem~\ref{lmm:yl-P1y} is a direct consequence of the following two results
      from number theory.
      \begin{lmm}{\normalfont \cite{yaglom1964challenging}}
      \label{lmm:yaglom}
      If $m > \binom{k+1}{2}$ then
      the number of partitions of $m$ in $k$ distinct parts
      equals the number of partitions of $m-\binom{k+1}{2}$ into
      at most $k$ parts (not necessarily distinct).
      \end{lmm}

      \begin{lmm}{\normalfont \cite{andrews1998theory,hardy1979introduction}}
      \label{lmm:partitionsk}
      Let $p_{\leq k}(m)$ be the number of partitions of $m$ in parts that do not exceed $k$.
      Then $p_{\leq k}(m)$ equals the number of partitions of $m$ into at most $k$ parts and
      \begin{equation*}
      %\label{eq:139}
      P_{\leq k}(q) = \sum_{m\geq0} p_{\leq k}(m) \, q^m =
      \prod_{i=1}^k \frac{1}{1-q^i} \,.
      \end{equation*}
      \end{lmm}

      % Lemma~\ref{lmm:partitionsk} is relatively easy to prove. Using Ferrers
      % diagrams, it is readily seen
      % that a partition in parts that do not exceed $k$ and a partition into
      % at most $k$ parts are conjugate
      % to one another;
      % see~\cite{andrews1998theory,hardy1979introduction} for more details.

      \begin{proof}[Proof of Theorem~\ref{lmm:yl-P1y}]
      Using Lemma~\ref{lmm:yaglom} and~\ref{lmm:partitionsk} we can recast~\eqref{eq:129} as
      \begin{align*}
      P(1,y) &= 1 + \sum_{k\geq1} \rho^k (y-1)^k \left[\sum_{m\geq{k+1 \choose 2}}
      d(m;k) q^m \right],\\
      &= 1 + \sum_{k\geq1} \rho^k (y-1)^k q^{{k+1 \choose 2}}\left[1+\sum_{m>0}
      d\left(m+{k+1 \choose 2};k\right) q^m \right],\\
      &= 1 + \sum_{k\geq1} \rho^k (y-1)^k q^{{k+1 \choose 2}}\left[1+\sum_{m>0}
      p_{\leq k}(m) q^m \right],\\
      &= 1 + \sum_{k\geq1} \rho^k (y-1)^k q^{{k+1 \choose 2}} \prod_{i=1}^k \frac{1}{1-q^i} \,,\\
      &= \sum_{k\geq0} \frac{\rho^k q^{{k+1 \choose 2}} (y-1)^k}{\prod_{i=1}^k [1-q^i] }\,,
      \end{align*}
      where, as usual, $p_{\leq k}(0)=1$ and $\prod_{i=a}^b f_k = 1$ when $b<a$.

      The second part of the Theorem is proved from~\eqref{eq:20} as follows:
      \begin{align*}
      \marginal &= \frac{1}{l!} \frac{d^l}{dy^l} P(1,y)\Big\vert_{y=0}\,,\\
      &= \sum_{k\geq l} \frac{\rho^k q^{{k+1 \choose 2}}}{\prod_{i=1}^k [1-q^i] }
      \frac{1}{l!} \frac{d^l}{dy^l} (y-1)^k \Big\vert_{y=0}\,,\\
      &= \sum_{k\geq l} \frac{(-1)^{k-l}\rho^k q^{{k+1 \choose 2}} {k\choose l}}{\prod_{i=1}^k [1-q^i] } \,.
      \end{align*}
      \end{proof}

      Now we use~\eqref{eq:138} to obtain an upper bound on $\marginal$:
      let $(q;q)_l = \prod_{i=1}^l \left[1-q^i\right]$
      be the $q$-Pochammer symbol of the pair $(q,q)$.
      The following inequalities hold:
      \begin{align*}
      {m+l \choose l} &= \prod_{k=1}^m \left( 1+\frac{l}{k} \right) \leq
      (1+l)^m \,,\\
      \prod_{i=l+1}^{l+m} \left[1-q^i\right] &\geq \prod_{i=1}^{m}
      \left[1-q^i\right] = (q;q)_m \,,\\
      (q;q)_l &\geq (q;q)_\infty\,.
      \end{align*}
      Thus,
      \begin{align}
      \marginal &= \frac{\rho^l\,q^{{l+1 \choose 2}}}{(q;q)_l}
      \sum_{m\geq0} \frac{(-\rho)^m q^{(l+1)m} q^{{m \choose 2}} \, {m+l \choose
      l}}{\prod_{i=l+1}^{l+m} [1-q^i]}\;, \nonumber\\%\label{eq:143}
      &\leq \frac{\rho^l\,q^{{l+1 \choose 2}}}{(q;q)_\infty}
      \sum_{m\geq0} \frac{ q^{{m \choose 2}} [\rho\, q^{l+1} \, (l+1)]^m}
      { (q;q)_m }\;, \label{eq:143a}\\
      &= \rho^l\,q^{{l+1 \choose 2}} \frac{ \prod_{k\geq0} \left[
      1 + q^{k+l+1} \rho (l+1) \right] }{(q;q)_\infty} \,, \label{eq:143}
      \end{align}
      where from~\eqref{eq:143a} to~\eqref{eq:143} we have used the
      properties of $q$-ascending factorials and $q$-binomial
      coefficients~\cite{gasper2004basic}.
      If $l$ is sufficiently large then $q^{l}\rho (l+1) \leq 1$,
      and~\eqref{eq:143} yields
      \begin{equation}
      \label{eq:24}
      \marginal \leq \rho^l\,q^{{l+1 \choose 2}}
      \frac{(-q;q)_\infty}{(q;q)_\infty} \,.
      \end{equation}
      \begin{rmk}
      Theorem~\ref{lmm:yl-P1y} and~\eqref{eq:24} show that, asymptotically
      in $l$, the leading order of \marginal{} is $\rho^l q^{{l+1 \choose
      2}}$. This fact can be directly implied from
      arrival process~\eqref{eq:eda1}. In fact, the most likely way to
      have $l$ late customers ($l$ large) is that each of the customers
      originally scheduled in the interval $[t-l, t)$ do not balk and are late at time $t$, an event of probability $\rho q^l \, \rho q^{l-1} \cdots \rho q =
      \rho^l q^{{l+1\choose 2}}$.
      \end{rmk}

      Since $P_{n,l} \leq \marginal$, we have just obtained the
      following asymptotic result:
      \begin{thm}
      \label{crl:1}
      Uniformly in $n$, the equilibrium distribution $P_{n,l}$ decays
      super-exponentially fast in $l$. More precisely,
      \begin{equation}
      \label{eq:21}
      P_{n,l} = O\left(\rho^l q^{{l+1 \choose 2}} \right)
      \quad \text{for } l \to \infty\,.
      \end{equation}
      \end{thm}

      As a matter of fact the super-exponential decay of $P_{n,l}$
      may be proved asymptotically in $n$.
      Let us consider the auxiliary process $\alpha_t = n_t +
      l_t$, which we have already encountered in the proof of
      Lemma~\ref{rmk:uniqueness}. There we have interpreted $\alpha_t$
      as the diagonal in the quarter plane where the point $(n_t,l_t)$ lies on.
      Under equilibrium conditions, the probability of finding the system on
      the $a$th diagonal is just
      \begin{equation*}
      p_a \equiv P(\alpha_t = a) = \sum_{\substack{n,l \geq 0\\n+l=a}} P_{n,l} \,,\quad a \geq 0 \,.
      \end{equation*}
      It is straightforward to prove that the generating function of $p_a$ is $P(z,z)$.
      % Then, we have the following:
      % \begin{lmm}
      % \label{lmm:6}
      % The generating function of $p_a$ is $P(z,z)$.
      % \end{lmm}
      % \begin{proof}
      % Let $A(z)=\sum_{a\geq0} p_a \, z^a$ be the generating function of $p_\alpha$.
      % Then,
      % \begin{equation*}
      % A(z) = \sum_{a\geq0} p_a \, z^a =
      % \sum_{a\geq0} \sum_{\substack{n,l \geq 0\\n+l=a}}
      % P_{n,l} \, z^{n+l} = \sum_{n,l\geq0} P_{n,l} \, z^n\,z^l= P(z,z) \,.
      % \end{equation*}
      % \end{proof}

      Substituting $y=z$ into~\eqref{eq:eda4} yields
      \begin{equation}
      \label{eq:22}
      P(z,z) = \frac{1+\rho(z-1)}{1-\rho} P(0,z) \,.
      \end{equation}
      \begin{rmk}
      Equation~\eqref{eq:22} gives an interesting connection between the
      equilibrium distribution of the quantity $\alpha_t = n_t+l_t$ and the
      stationary probability of having $l$ late
      customers given that the queue is void.
      Figure~\ref{fig:quarterplane} shows that the latter event drives the
      dynamic of $\alpha_t$ through an independent Bernoulli random
      variable with parameter $\rho$, which explains the factor $1+\rho(z-1)$.

      \end{rmk}

      From~\eqref{eq:22} we can compute as follows $p_a$ in terms of $P_{0,a}$:
      %\begin{fleqn}
      \begin{align}
      p_a &= \frac{1}{a!}\,\frac{d^a}{d z^a}P(z,z)
      \Big\vert_{z=0} \nonumber \\
      &= \frac{1}{a!}\,\Big[\frac{1+\rho(z-1)}{1-\rho}\frac{d^a}{d
      z^a}P(0,z)
      + \frac{\rho\,a}{1-\rho} \frac{d^{a-1}}{d z^{a-1}}P(0,z)
      \Big]_{z=0} \,, \nonumber \\
      \label{eq:132}
      &= P_{0,a} + \frac{\rho}{1-\rho}P_{0,a-1}\,.
      \end{align}
      For $a=n+l$, formulas~\eqref{eq:21} and~\eqref{eq:132} then yield
      \begin{equation}
      \label{eq:15}
      P_{n,l} \leq p_a = P_{0,a} +
      \frac{\rho}{1-\rho}P_{0,a-1} = O(\rho^{a} q^{\binom{a}{2}}) \,.
      \end{equation}
      Therefore, the following asymptotic result holds:
      \begin{thm}
      The equilibrium distribution $P_{n,l}$ decays super-exponentially
      fast as either $n\to\infty$ or $l\to\infty$. More precisely,
      \begin{equation}
      \label{eq:super-exp-decay}
      P_{n,l} = O\left(\rho^{n+l} q^{{n+l \choose 2}} \right)
      \quad \text{for } n,l \to \infty\,.
      \end{equation}
      \end{thm}

    \section{Numerical approximation of the joint
    stationary measure}
    \label{sec:numerics}

      In this final section we examine the possibility to approximately
      compute the joint stationary distribution $P_{n,l}$. Due to the very
      broad range of applications of the queueing model $EDA/D/1$, an
      efficient approximate computation of the solution may prove itself
      crucial in contexts where practical solutions are needed.

      In Section~\ref{sec:marg-distr-late} we have shown that the joint
      stationary probability $P_{n,l}$ decreases super-exponentially fast in
      the limit of either $n, l \to \infty$. The natural question arising is
      then whether a bare truncation of the infinite
      linear system~\eqref{eq:eda95}-\eqref{eq:eda97} is sufficient to obtain a
      satisfactory numerical expression of $P_{n,l}$.
      As we will see, in this case the answer is positive due
      to~\eqref{eq:super-exp-decay}.

      \smallskip

      We truncate the infinite system of balance equations~\eqref{eq:eda95}-\eqref{eq:eda97} by fixing an integer $\alpha_{\max}$ and imposing $P_{n,l} = 0$ for $n+l > \alpha_{\max}$. For the purpose of simplifying the notation, we map the quarter plane $\{n,l \in \mathbb{N}\times\mathbb{N}\}$ onto the set of non-negative integers.
      This way, we can relabel the unknowns $P_{n,l}$ as $\pi_i$ and recast~\eqref{eq:eda95}-\eqref{eq:eda97} as $\pi = \pi Q$. For the details of both mapping and relabeling, see Appendix~\ref{sec:map}.
      Next, we map $\{n,l \in \mathbb{N}\times\mathbb{N} \text{ such that } n+l \leq \alpha_{\max}\}$ onto the set of non-negative integers $\{0,1,\ldots,k_{\max}\}$, where $k_{\max} = {\alpha_{\max} + 1 \choose 2}$.
      We want to consider the truncated system
      \begin{equation}
        \begin{cases}
          & \tilde{\pi}_i = \sum_{j=0}^{k_{\max}} \tilde{\pi}_j Q_{j,i} \qquad i = 1,2,\ldots,k_{\max}-1 \,,\\
          & \sum_{j=0}^{k_{\max}} \tilde{\pi}_j = 1 \,.
        \end{cases}
      \end{equation}

      The idea we present is not new and has been already discussed, for
      instance in~\cite{tijms1994stochastic} for stationary distributions
      with geometric tail.
      As shown in Appendix~\ref{sec:truncated_system}, there exists a sequence $\{\varepsilon_j\}$ such that
      \begin{align*}
      & \sum_{i > k_{\max}} \pi_i \, Q_{i,j} \leq \varepsilon_j \,, \qquad
      j=0,1,\ldots,k_{\max}-1\,,\\
      & \sum_{i > k_{\max}} \pi_i \leq \varepsilon_{k_{\max}} \,.
      \end{align*}
      The following a priori estimate of the error introduced by the truncation can be obtained from \emph{perturbation theory}~\cite[\S 2.6.2]{golub2012matrix}:
      \begin{equation}
      \label{eq:38}
      \sum_{j=0}^{k_{\max}} |\pi_j - \tilde{\pi}_j| \leq \kappa(A) \sum_{i=0}^k
      \varepsilon_i \,,
      \end{equation}
      where
      $A$ is the $k_{\max} \times k_{\max}$ matrix
      \begin{equation}
      \label{eq:34}
      A_{ij} =
      \begin{cases}
      \delta_{i,j} - Q_{i,j} \,, &\quad
      i=0,1,\ldots,k_{\max}-1\,,\\
      1 & \quad i=k_{\max} \,,
      \end{cases}
      \end{equation}
      $\delta_{i,j}$ is the usual Kronecker's delta, and
      \[\kappa(A) = \Vert A \Vert_1 \Vert A^{-1} \Vert_1 \]
      is the norm-1 condition number of the matrix $A$.

\begin{figure}[tbp]
      \centering
      \includegraphics[width=.85\textwidth]{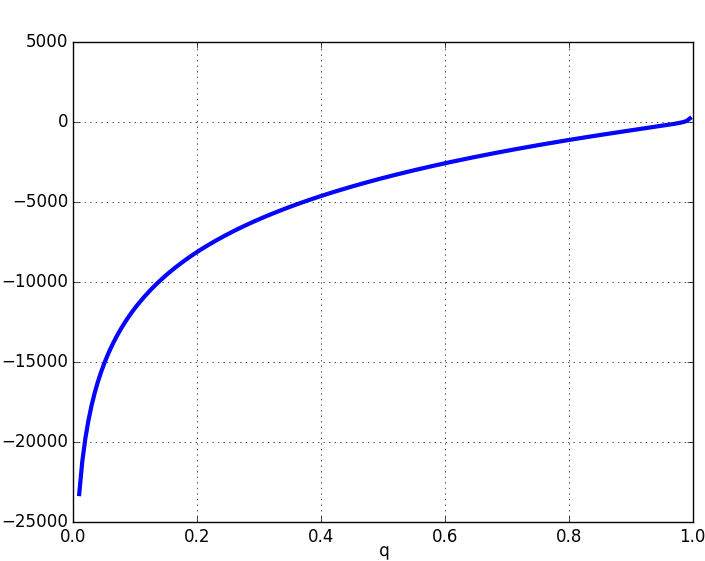}
      \caption{Behaviour of $\log\left(\frac{(-q;q)_\infty}{(q;q)_\infty}
      q^{{\alphamax \choose 2}}\right)$ as a function of $q$ for $\alpha_{\max}=100$.}
      \label{fig:logbound}
      \end{figure}

      From Appendix~\ref{sec:truncated_system},
      \begin{equation}
      \label{eq:36}
      \sum_{j=0}^{\kmax} \varepsilon_j \leq 2 \alpha_{\max} \frac{(-q;q)_\infty}{(q;q)_\infty}
      \rho^{\alphamax+1} q^{{\alphamax + 1 \choose 2}} \leq 2
      \alpha_{\max} \frac{(-q;q)_\infty}{(q;q)_\infty} q^{{\alphamax + 1 \choose 2}}\,.
      \end{equation}
      Figure~\ref{fig:logbound} shows the behaviour of
      $\log\left(\frac{(-q;q)_\infty}{(q;q)_\infty} q^{{\alphamax \choose 2}}\right)$ as a function of $q$ for
      $\alpha_{\max}=100$.
      Therefore, unless the condition number of the
      matrix $A$ is very large, we expect that a truncation at the level
      $\alpha_{\max}=100$ will give a very good approximation of the
      stationary probabilities of $EDA/D/1$.

\begin{figure}[tbp]
      \centering
      \includegraphics[width=.8\textwidth]{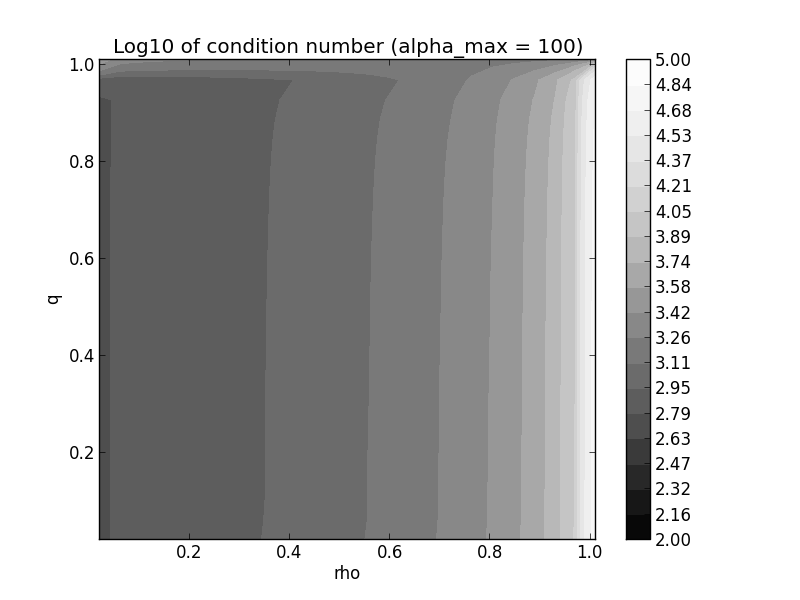}
      \caption{$\log_{10}\kappa(A)$ for $\alpha_{\max} = 100$ when
      $\rho$ and $q$
      vary between $0.0$ and $0.99$ in $25$ steps.
      The condition number is of order
      $10^5$ at most in this region of the parameters.}
      \label{fig:CondNo_alpha100}
      \end{figure}

      Estimating the condition number of a
      matrix is a notably difficult problem and a
      vast literature exists on this topic, see
      \eg~\cite{golub1965calculating,qi1984some, johnson1989gersgorin}.
      Since the aim of the present section is showing that a bare
      truncation of the balance
      equations~\eqref{eq:eda120}--\eqref{eq:eda121} may be sufficient for
      an approximate computation of $P_{n,l}$, we fall back on
      numerical computations.
      Figure~\ref{fig:CondNo_alpha100} displays the value of
      $\kappa(A)$ in the $\rho q$-plane when $\alpha_{\max} = 100$.
      We see that the condition number is not larger than $10^5$ for
      $\rho, q\leq 0.99$.

      Table~\ref{tab:bound} gives the numerical value of the right-hand side
      of~\eqref{eq:36} for $q$ between $0.9$ and $0.99$, and
      $\alpha_{\max}=100$. Comparison of Table~\ref{tab:bound} with
      Figure~\ref{fig:CondNo_alpha100} shows that, uniformly in $\rho\leq0.99$, the
      a priori norm-1 approximation error is less than $10^{-12}$ for $q$ up
      to $0.98$.
      \begin{table}[tbp]
      \centering
      \caption{Value of $2\alpha_{\max}\frac{(-q;q)_\infty}{(q;q)_\infty}
      q^{{\alphamax \choose 2}}$ for $0.90 \leq q \leq 0.99$ and
      $\alpha_{\max}=100$.}
      \begin{tabular}{*{10}{c}}
      \toprule
      $0.90$ & $0.91$ & $0.92$ & $0.93$ & $0.94$\\
      $1.4\times10^{-224}$ & $9.8\times10^{-200}$ & $5.1\times10^{-175}$
       & $2.4\times10^{-150}$ &  $1.3\times10^{-125}$\\
      \midrule
      $0.95$ & $0.96$ & $0.97$ & $0.98$ & $0.99$\\
      $1.3\times10^{-100}$ & $5.6\times10^{-75}$ & $8.7\times10^{-48}$
       & $2.3\times10^{-16}$ & $2.6\times10^{33}$\\
      \bottomrule
      \end{tabular}
      \label{tab:bound}
      \end{table}
      \begin{rmk}
      For air traffic applications, $q=0.98$ corresponds to typical delays
      of the order of one hour. The same value of $q$ for other transport
      systems, \eg{} trains or buses, correspond to even higher
      delays. Consider also that typical values of $\rho$ in extremely
      congested systems, \eg{} London Heathrow Airport, do not exceed
      $0.98$, see~\cite{cills}.
      Therefore, the approximation scheme presented in this
      section is very fit for real life applications.
      \end{rmk}

\begin{figure}[tbp]
      \centering
      \includegraphics[width=.78\textwidth]{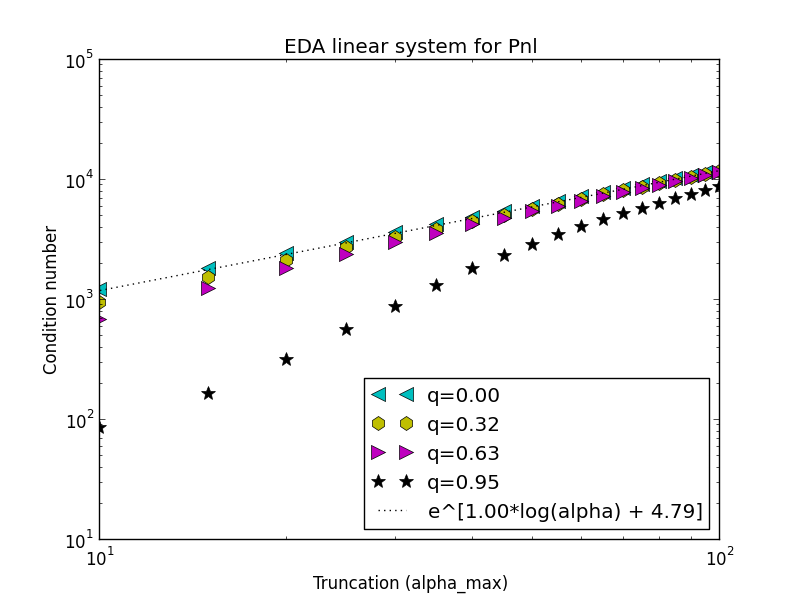}
      \caption{Log-log plot of the condition number
      $\kappa(A)$ for $\alpha_{\max}$ varying between $10$ and
      $100$, $\rho=0.95$, and different values of $q$. The
      curve slope for $q=0.0$ is $1.0$.}
      \label{fig:CondNo_loglog_alpha100}
      \end{figure}

      Finally, Figure~\ref{fig:CondNo_alpha100} suggests that the system
      condition number is decreasing in $q$ for fixed $\rho$.
      Figure~\ref{fig:CondNo_loglog_alpha100} validates this insight by showing a
      log-log plot of the condition number for $\alpha_{\max}$ between $10$ and
      $100$ and $\rho = 0.95$.
      The very same figure also suggests that the condition number of
      the system may grow polynomially with $\alpha_{\max}$. In particular,
      $\kappa(A)$ seems to grow linearly in
      $\alpha_{\max}$ for $\rho =
      0.95$ and $q = 0.0$.

      \section{Conclusions}
      \label{sec:conclusions}

      In this paper we have addressed a single-server queueing system
      with deterministic service time and \emph{exponentially delayed arrivals}.
      The point process describing these arrivals dates back to the
      '50s of the past century and was studied by Kendall and others.

      We have described the model as a bivariate Markov chain, proved that
      the latter is ergodic,
      wrote the balance equations of the stationary distribution,
      and found a functional equation for the bivariate
      generating function.
      Then we have focused on the marginal distribution
      of the number of late customers and found its exact expression.
      This intermediate step has enabled the fundamental
      result on the super-exponential decay of the joint stationary distribution.
      The characterisation of the asymptotic behaviour
      has finally led us to show that the solution to
      the balance equations can be approximately computed in
      a simple yet very accurate way.

      In spite of the big efforts we have put forward to find the solution
      to the functional equation~\eqref{eq:eda4},
      the complete solution of the problem is still out of reach.
      An expansion in powers of $q$
      seems to be a promising approach to obtain (at least) an approximate
      expression of the bivariate generating function.
      This method allows to set up a recursive
      scheme to compute the coefficients of the power series,
      see~\cite{oldEDA:2013,lancia13looking}.
      Unfortunately, the computations are quite involved and need
      some additional work to be refined. This will be the subject
      of further research and the topic of an upcoming paper.

      Figures~\ref{fig:logbound}--\ref{fig:spyAbar} were obtained using
      Python~2.7.9,
      numpy~1.9.1,
      scipy~0.15.1
      matplotlib~1.4.2, and
      mpmath~0.19.
      The code to generate them is freely available on GitHub
      at the following address: \url{https://github.com/clancia/EDA}.

      \appendix

      \section{Map of the quarter plane onto non-negative integers} % (fold)
      \label{sec:map}

        Define the map $F$ of the quarter plane onto
        the set of the non-negative integers and its inverse $G$:
        \begin{align*}
        & F(n,l) \mapsto m= {n+l \choose 2} + l \,,\\
        & g(m) = \max\left\{j \geq 0 \text{ such that } {j+1
        \choose 2} \leq m\right\}  = \floor*{\frac{-1+\sqrt{1+8m}}{2}} \,,\\
        & G(m) = \left(g(m), m-g(m) \right) \,,
        \end{align*}
        where $\floor*{\cdot}$ denotes the lower integer part, \ie{} the
        \emph{floor} operation.
        Fixed a positive integer $\alpha_{\max}$, let $k_{\max} =
        {\alpha_{\max} + 1 \choose 2}$.
        Define the $k_{\max} \times k_{\max}$ matrix
        \begin{equation}
        \label{eq:33}
        A_{ij} =
        \begin{cases}
        \delta_{i,j} - \mathcal{P}(G(i), G(j))\,, &\quad
        i=0,1,\ldots,k_{\max}-1\,,\\
        1 & \quad i=k_{\max}\,,
        \end{cases}
        \end{equation}
        where $\mathcal{P}(\cdot,\cdot)$ is defined by~\eqref{eq:eda120}-\eqref{eq:eda121}, and
        \begin{equation}
        \label{eq:37}
        \pi_j = P_{g(j), j-g(j)} \,, \qquad j=0,1,\ldots,k_{\max} \,.
        \end{equation}

\begin{figure}[tbp]
        \centering
        \includegraphics[width=.65\textwidth]{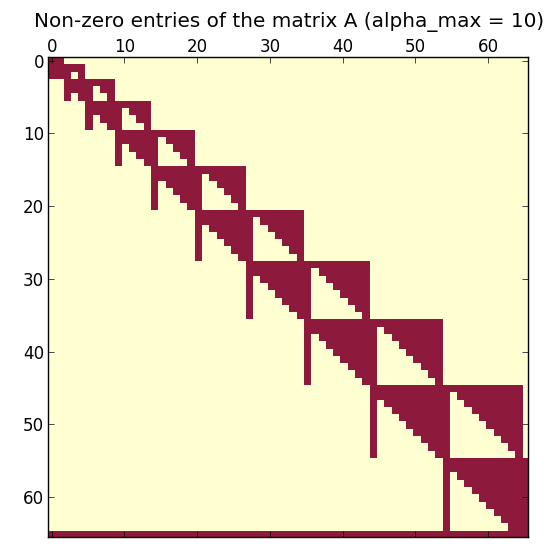}
        \caption{Sparsity structure of the matrix $A$,
        non-zero elements in dark colour ($\alpha_{\max}$ is set to $10$
        for readability).}
        \label{fig:spyAbar}
        \end{figure}

        \begin{rmk}
        The matrix~\eqref{eq:33} is rather sparse, as shown by
        Figure~\ref{fig:spyAbar}. We recommend to exploit this
        property by using sparse storage formats and dedicated libraries
        when $q$ is set larger than $0.98$.
        In this regime $\alpha_{\max}=100$ could be no longer
        sufficient to achieve a good approximation of $P_{n,l}$, but
        enlarging $\alpha_{\max}$ while using dense formats could quickly
        lead to memory shortage and a severe computational slowdown.
        \end{rmk}

      % section map (end)

      \section{Truncated System} % (fold)
      \label{sec:truncated_system}

      By product of~\eqref{eq:super-exp-decay} and direct inspection
      of~\eqref{eq:eda120}--\eqref{eq:eda121}, for $i+j = \alpha_{\max}$,
      \begin{equation}
      \label{eq:39}
      \sum_{\substack{n,l\geq0\\n+l>\alphamax}} P_{n,l} \, \mathcal{P}((n,l),
      (i,j))  \leq (1-\rho) \frac{(-q;q)_\infty}{(q;q)_\infty}
      \rho^{\alphamax+1} q^{{\alphamax + 1 \choose 2}} \,,
      \end{equation}
      while for $i+j < \alpha_{\max}$,
      \begin{equation}
      \label{eq:40}
      \sum_{\substack{n,l\geq0\\n+l>\alphamax}} P_{n,l} \, \mathcal{P}((n,l),
      (i,j))  = 0 \,,
      \end{equation}
      where $P_{n,l}$ and $\mathcal{P}(\cdot, \cdot)$ are defined
      by~\eqref{eq:eda120}--\eqref{eq:eda121} and~\eqref{eq:pnldef},
      respectively. Also,
      \begin{equation}
      \label{eq:41}
      \sum_{\substack{n,l\geq0\\n+l>\alphamax}} P_{n,l}  \leq 2 \frac{(-q;q)_\infty}{(q;q)_\infty}
      \rho^{\alphamax+1} q^{{\alphamax + 1 \choose 2}}\,.
      \end{equation}

      % section truncated_system (end)

    \section*{Acknowledgements}
      B.S. and C.L. have been supported by PRIN 2012 ``Problemi matematici in teoria cinetica ed applicazioni''.
      C.L.\ thanks the Mathematical Institute of Leiden University for the warm hospitality.
      G.G. and S.N.\ appreciate
      the Mathematics Department of the University of
      Rome `Tor Vergata' for the kind support.

    % bib stuff
    %\nocite{*}
    \addtocontents{toc}{\protect\vspace{\beforebibskip}}
    \addcontentsline{toc}{section}{\refname}
    \bibliographystyle{siam}
    \bibliography{eda}
\end{document}